\documentclass[a4paper]{article}
\usepackage{amsmath,amssymb}
\usepackage[latin1]{inputenc}
\usepackage[dvips]{epsfig}

\newcommand{\ZZ}{\mathbb{Z}}
\newcommand{\CC}{\mathbb{C}}
\newcommand{\For}{\mathsf{F}}
\newcommand{\mon}{\mathsf{M}}
\newcommand{\for}{\mathsf{for}}
\newcommand{\souche}{\operatorname{Stump}}
\newcommand{\crk}{\operatorname{crk}}
\newcommand{\tree}{\mathsf{Tree}}
\newcommand{\zero}{\widehat{0}}
\newcommand{\depth}{\operatorname{Depth}}

\newcommand{\der}{\operatorname{Der}}
\newcommand{\dq}{\mathsf{Q}}

\newtheorem{prop}{Proposition}
\newtheorem{theo}{Theorem}

\newtheorem{lemma}{Lemma}

\newenvironment{proof}{\begin{trivlist}\item{\bf{Proof.}}}
  {\hfill\rule{2mm}{2mm}\end{trivlist}}

\title{Free hyperplane arrangements associated to labeled rooted trees}

\date{\today}

\author{Frédéric Chapoton}

\begin{document}

\maketitle

\begin{abstract}
  To each labeled rooted tree is associated a hyperplane arrangement,
  which is free with exponents given by the depths of the vertices of
  this tree. The intersection lattices of these arrangements are
  described through posets of forests. These posets are used to define
  coalgebras, whose dual algebras are shown to have a simple
  presentation by generators and relations.
\end{abstract}

\setcounter{section}{-1}

\section{Introduction}

This article is centered on labeled rooted trees. This kind of tree is
very classical in combinatorics, since the enumerative results of
Cayley \cite{cayley89}. More recently, it has surfaced, with more
algebraic structure, in the study of pre-Lie algebras \cite{rooted}
and in relation with a Hopf algebra in renormalization
\cite{conneskreimer}.

Here, we deal with other algebraic and geometric aspects of rooted
trees. It is not yet clear if these new aspects are related to the
previous ones.

The starting point is the definition of a hyperplane arrangement for
each individual rooted tree. One can consider a rooted tree as a
poset, and from this viewpoint directly follow the equations of the
arrangement. Apart from being quite simple, these arrangements have
the remarkable property of being free, in the sense of Saito
\cite{saito}. More precisely, the arrangement associated to a rooted
tree is free with exponents given by the depths of the vertices of
this tree.

After the first version of this article was completed, I learned from
R. Stanley that this result is a consequence of the theory of
graphical arrangements \cite[\S 3, Th. 3.3]{edelrein} and
supersolvable lattices \cite[Prop.  2.8, Ex. 4.6]{stanley72}. A
graphical arrangement is known to be free if and only if the
corresponding graph is chordal. As the comparability graphs for rooted
trees are indeed chordal, one can recover in this way Theorem
\ref{liberte}. My proof of freeness avoids using these general
theories and gives explicit information on logarithmic vector fields
and differential forms.

Besides, it turns out that the intersection lattices of these
hyperplane arrangements admit a neat combinatorial description in
terms of forests of labeled rooted trees. This leads to the definition
of a partial order on the set of forests labeled by a finite set,
which contains intersection lattices as intervals.

Next, inspired by similarity with the study of binary leaf-labeled
trees made in \cite{bessel,besselposet}, the posets of forests are
used to define coalgebras on forests. The dual algebras of these
coalgebras are then shown to have a simple presentation by generators
and relations.

\smallskip

Thanks to David Bessis for helpful discussions on the topology of the
complement.

\section{The category of labeled rooted trees}

\label{category}

Let $I$ be a finite set. A \textit{tree} on $I$ is a connected, simply
connected graph on the vertex set $I$ endowed with a distinguished
vertex called the root.

It will be convenient to use the following equivalent definition. A
tree on $I$ is given by a partial order relation on the finite set $I$
with a unique minimal element (the root) and such that the interval
between the root and any element of $I$ is a chain. The order relation
corresponding to a tree $T$ is denoted by $\leq_T$.

\begin{figure}
  \begin{center}
    \leavevmode 
    \epsfig{file=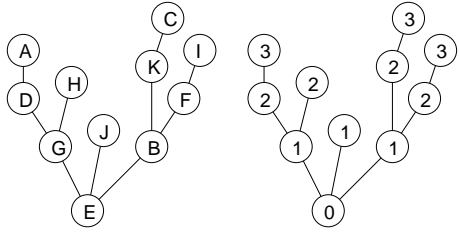,width=6cm} 
    \caption{A labeled rooted tree and the depths of its vertices.}
    \label{fig-prof}
  \end{center}
\end{figure}

The \textit{depth} of a vertex $i$ is defined to be the number of
edges in the maximal chain from the root to $i$ and is denoted by
$\depth(i)$, see Fig. \ref{fig-prof}.

By convention, trees are drawn with their root at the bottom and edges
are oriented towards the root, in the decreasing direction for the
poset structure. The \textit{valence} of a vertex is the number of its
incoming edges.

A \textit{linear tree} is a tree whose vertices have valence at most
$1$.

The category $\tree$ is defined as follows. Its objects are the trees
on $I$ for all finite sets $I$. Morphisms from a tree on $I$ to a tree
on $J$ are the maps from $I$ to $J$ which are morphisms of posets.

\section{The hyperplane arrangement of a rooted tree}

Let $I$ be a finite set. Consider the vector space $\CC ^I$ with
coordinates $x_i$. Let $T$ be a tree on $I$. The equations
\begin{equation}
  x_i=x_j \quad \text{if} \quad i \leq_T j
\end{equation}
define a central hyperplane arrangement in $\CC ^I$ denoted by $H_T$.

Note that the complement of the union of these hyperplanes is an
algebraic variety defined over $\ZZ$ which depends functorially from
the tree in the category $\tree$ of rooted labeled trees defined
above.

Let $\dq_T$ be the product of all the equations $x_i-x_j$ for $i
\leq_T j$. It is called the \textit{defining form} of the arrangement $H_T$.

Here is the first main result of this article.

\begin{theo}
  \label{liberte}
  The arrangement $H_T$ is free. Its exponents are the depths of the
  vertices of $T$.
\end{theo}

Remark that the arrangement $H_T$ for $T$ a linear tree on $I$ is
simply the braid arrangement given by the reflection hyperplanes of
the symmetric group of $I$. Therefore the classical freeness of this
braid arrangement is recovered as a special case of Theorem \ref{liberte}.

The proof of Theorem \ref{liberte} uses the Saito criterion of
freeness \cite[Th. 1.8 (ii)]{saito} involving logarithmic vector
fields and is the object of the next section.

\section{Logarithmic vector fields}

First one can associate a vector field $\theta_i$ to each vertex $i$
of $T$. Define the vector field $\theta_i$ on $\CC^I$ by
\begin{equation}
  \theta_i =\sum_{j \geq_T i}\left(\prod_{k <_T i}(x_k-x_j)\right)\partial_j.
\end{equation}
The degree of $\theta_i$ is one less than the depth of $i$ in $T$.

This definition can be restated as
\begin{equation}
  \label{thetatriangle}
  \theta_i(x_j)=
  \begin{cases}
    \prod_{k <_T i}(x_k-x_j) &\text{ if }i\leq_T j,\\
    0 &\text{ otherwise.}
  \end{cases}
\end{equation} 

\smallskip

Let $\der(H_T)$ be the space of logarithmic vector fields for the
arrangement $H_T$, \textit{i.e.} the space of polynomial vector fields
$\theta$ on $\CC ^I$ such that $\alpha$ divides $\theta(\alpha)$ for
all hyperplanes $\alpha$ of $H_T$.

\begin{prop}
  The vector fields $\theta_i$ belong to $\der(H_T)$. 
\end{prop}

\begin{proof}
  Let $j$ and $k$ be distinct vertices of $T$ with $j \leq_T k$.
  
  Assume first that neither $i \leq_T j$ nor $i \leq_T k$. Then
  $\theta_i(x_j-x_k)$ is zero.
  
  Assume then that $i \not\leq_T j$ but $i \leq_T k$. As $T$ is a
  tree, these conditions together with $j \leq_T k$ implies that
  $j<_T i$. It follows that $\theta_i(x_j-x_k)=-\prod_{\ell <_T
    i}(x_\ell-x_k)$ is divisible by $x_j-x_k$.
  
  Assume now that $i \leq_T j$, then also $i \leq_T k$ by
  transitivity. Then $\theta_i(x_j-x_k)=\prod_{\ell <_T
    i}(x_\ell-x_j)-\prod_{\ell <_T i}(x_\ell-x_k)$ is divisible by
  $x_j-x_k$, because this expression vanishes when $x_j=x_k$.

  Therefore $\theta_i(x_j-x_k)$ is always divisible by $x_j-x_k$.
\end{proof}

One can now prove Theorem \ref{liberte}.

\begin{proof}
  By Formula (\ref{thetatriangle}), the matrix
  $\Theta=[\theta_i(x_j)]$ is triangular for the poset structure on
  $I$ given by $\leq_T$. The determinant of the matrix $\Theta$ is
  therefore
  \begin{equation}
    \prod_{(k,i) \atop {k<_T i}}(x_k-x_i),
  \end{equation}
  which is the defining form $\dq_T$ of the arrangement $H_T$. Hence
  the Saito criterion applies and the arrangement $H_T$ is free. The
  exponents are by definition one more than the degrees of the vector
  fields $\theta_i$, which are one less than the depths of the
  vertices of $T$.
\end{proof}

Freeness of the arrangement $H_T$ means that $\der(H_T)$ is a free
module over the ring of polynomials on $\CC^I$. The proof via the
Saito criterion implies further that the set of vector fields
$(\theta_i)_{i \in I}$ is a basis of $\der(H_T)$.

\section{Logarithmic differential forms}

Define a differential $1$-form $\omega_i$ for each vertex $i$ of $T$
by
\begin{equation}
  \omega_i=
  \sum_{j \leq_T i} \left( \prod_{{k\leq_T i}\atop{k\not=j}} 
    \frac{1}{x_k-x_j}\right) dx_j.
\end{equation}

Recall the definition of the space of logarithmic differential forms
for an hyperplane arrangement \cite{saito}. A differential form
$\omega$ with coefficients in the field of rational functions on
$\CC^I$ is called a logarithmic differential form for $H_T$ if $\dq_T
\omega$ and $\dq_T d\omega$ are polynomial forms on $\CC ^I$.

\begin{prop}
  The forms $\omega_i$ are logarithmic $1$-forms for $H_T$.
\end{prop}
\begin{proof}
  It is clear that $\dq_T \omega_i$ is polynomial. On the other hand,
  one has
  \begin{equation}
   -d \omega_i=\sum_{j \leq_T i} 
    \sum_{k\leq_T i \atop {k\not=j}}\left( \prod_{{\ell\leq_T i}\atop{{\ell\not=j}\atop{\ell\not=k}}} \frac{1}{x_\ell-x_j}\right)
    \frac{1}{(x_k-x_j)^2} dx_k \wedge dx_j.
  \end{equation}
  Consider now the coefficient of $dx_k \wedge dx_j$ in $d \omega_i$.
  It is sufficient to prove that it has only a simple pole at
  $x_k=x_j$. This follows from the vanishing of
  \begin{equation}
    \left( \prod_{{\ell\leq_T i}\atop{{\ell\not=j}\atop{\ell\not=k}}} \frac{1}{x_\ell-x_j}\right)
    -\left( \prod_{{\ell\leq_T i}\atop{{\ell\not=j}\atop{\ell\not=k}}} \frac{1}{x_\ell-x_k}\right)
  \end{equation}
  when $x_k=x_j$.
\end{proof}

As the arrangement $H_T$ is free, it is known that the space of
logarithmic $1$-forms is a free module over the the ring of
polynomials on $\CC^I$. Furthermore this module is in duality with the
module of logarithmic vector fields by restriction of the duality
between vector fields and $1$-forms \cite{saito}.

More precisely, one has the following proposition in the case of $H_T$.

\begin{prop}
  \label{dualbasis}
  The set of logarithmic $1$-forms $(\omega_i)_{i \in I}$ is the dual
  basis to the basis $(\theta_i)_{i \in I}$.
\end{prop}
\begin{proof}
  One has 
  \begin{multline}
    \langle \omega_i, \theta_{i'} \rangle=\sum_{i'\leq_T j \leq_T
      i}\left( \prod_{ k\leq_T i \atop{k\not=j}} 
      \frac{1}{x_k-x_j}\right)\left( \prod_{j' <_T i'} 
      {x_{j'}-x_j}\right)\\
    =    \sum_{i'\leq_T j \leq_T i}\left( \prod_{i'\leq_T k\leq_T i
      \atop {k \not=j}} 
      \frac{1}{x_k-x_j}\right).
  \end{multline}
  The proposition now follows from Lemma \ref{chaindirac} below.
\end{proof}
\begin{lemma}
  \label{chaindirac}
  The sum 
  \begin{equation}
    \sum_{i'\leq_T j \leq_T i}\left( \prod_{i'\leq_T k \leq_T i \atop{k\not=j}} 
      \frac{1}{x_k-x_j}\right)
  \end{equation}
  equals $1$ if $i'=i$ and vanishes otherwise.
\end{lemma}
\begin{proof}
  This sum is clearly equal to $1$ when $i=i'$ and $0$ when $i' >_T
  i$. Assume now that $i'<_T i$. Then the sum is homogeneous of
  negative degree. As poles are at most of order one, it is sufficient
  to prove that all residues vanish. The residue at $x_{j_0}=x_{k_0}$
  is given by the value of
  \begin{equation}
    \left( \prod_{i'\leq_T k\leq_T i \atop{k\not=j_0 \atop{k\not=k_0}}} 
      \frac{1}{x_k-x_{j_0}}\right)-
    \left( \prod_{i'\leq_T k\leq_T i \atop{k\not=j_0 \atop{k\not=k_0}}} 
      \frac{1}{x_k-x_{k_0}}\right)
  \end{equation}
  at $x_{j_0}=x_{k_0}$, which is zero.
\end{proof}

\section{Topology of the complement}

\label{topologie}

Let $M_T$ denote the complexified complement of the union of all
hyperplanes of $H_T$. Using the general theory of supersolvable
arrangements \cite{stanley72,teraoss} and the results of
\cite{edelrein} on graphical arrangements, the complement $M_T$ can be
precisely described.

One can map a tree $T$ to the comparability graph $\Gamma_T$ for the
relation $\leq_T$. In this way, a rooted tree is seen as a special
case of a graph on $I$. Then $H_T$ is the graphical arrangement
corresponding to $\Gamma_T$. It is also easy to see that the graph
$\Gamma_T$ is chordal (see \cite{edelrein} for the definition).

Recall from \cite[\S 3]{edelrein} the notions of simplicial vertex and
vertex elimination ordering of a graph. For the graph $\Gamma_T$, one
has
\begin{lemma}
  The leaves of $T$ (maximal vertices for $\leq_T$) are simplicial
  vertices. Any total order extending the partial order $\leq_T$ on
  $I$ gives a vertex elimination ordering.
\end{lemma}
Such an ordering will be called a \textit{leaf-removal ordering}.

As $\Gamma_T$ is chordal, Theorem 3.3 of \cite{edelrein} implies that
$H_T$ is a supersolvable arrangement and that each leaf-removal
ordering gives a modular chain in $L_T$. By results of Terao
\cite{teraoss}, each modular chain in a supersolvable arrangement
gives rise to a description of the complement as an iterated fibration
of punctured $\CC$.

Hence the complement $M_T$ is an iterated fibration of the spaces $\CC
\setminus S_i$ for $i\in I$, where $S_i$ is a finite set of distinct
points in $\CC$ of cardinality the depth of $i$ in $T$. The order of
the fibration tower is given by the chosen leaf-removal ordering of
$I$. From this, one can deduce the following proposition.

\begin{prop}
  The space $M_T$ is a $K(\pi,1)$ space and its fundamental group is
  an iterated extension of free groups on $\depth(i)$ generators for
  $i\in I$, in the order given by any leaf-removal ordering of $I$.
\end{prop}

\section{Cohomology of the complement}

One can recall the results of Orlik and Solomon \cite{orlik} on the
cohomology of the complexified complement.

\begin{prop}
  The cohomology of the complex complement of $H_T$ is generated by
  the differential forms $d(\log(x_i-x_j))$ for $i \leq _T j$.
\end{prop}

In the same way, the similar results of Gelfand and Varchenko
\cite{gelfvarch} on the filtered algebra of locally constant integral
functions on the real complement of a real hyperplane arrangement can
be applied.

It would be interesting to find the relations for these algebras,
\textit{i.e.} to describe the dependent sets of hyperplanes. As a
first step in this direction, one can give a generating set of
relations between the linear forms $\alpha_{i,j}=x_i -x_j$ for $i
\leq_T j$.

\begin{lemma}
  Any linear relation between the linear forms $\alpha_{i,j}$ for $i
  \leq_T j$ is a linear combination of relations
  \begin{equation}
    \label{elementaire}
    \alpha_{i,j}+\alpha_{j,k}-\alpha_{i,k}=0,
  \end{equation}
  with $i \leq_T j \leq_T k$.
\end{lemma}

\begin{proof}
  Let $\ell$ be the root of $T$. Consider an arbitrary linear relation
  \begin{equation}
    \sum_{i \leq_T j} \lambda_{i,j} \alpha_{i,j}=0.
  \end{equation}
  This relation can be rewritten, using all relations
  \begin{equation}
    \label{relaracine}
    \alpha_{i,j}=\alpha_{\ell,j}-\alpha_{\ell,i},
  \end{equation}
  as a relation of the following kind:
  \begin{equation}
    \sum_{i \not= \ell} \mu_i \alpha_{\ell,i}=0.
  \end{equation}
  Now the coefficients $\mu_i$ vanish because the forms
  $\alpha_{\ell,i}$ are obviously linearly independent. Note that
  among relations (\ref{elementaire}) only relations
  (\ref{relaracine}) involving the root are really used.
\end{proof}

%%%%%%%%%%%%%%%%%%%%%%%%%%%%%%%%%%%%%%%%%%%%%%%%%%%%%%%%%%%%%%%%%%%%%%
                                % partie II 
%%%%%%%%%%%%%%%%%%%%%%%%%%%%%%%%%%%%%%%%%%%%%%%%%%%%%%%%%%%%%%%%%%%%%%

\section{The lattice of a labeled rooted tree}

Let $L_T$ be the intersection lattice of the arrangement $H_T$,
\textit{i.e.} the lattice of intersections of the hyperplanes of
$H_T$ for the reverse inclusion order.

\subsection{Characteristic polynomial}

The characteristic polynomial of $H_T$ can be deduced from a theorem
of Terao on free arrangements (Main Theorem of \cite{terao81}, see
also \cite[Th. 5.1]{sagan}). Recall that the characteristic polynomial
of the hyperplane arrangement $H_T$ is
\begin{equation}
  \chi_T(y)=\sum_{a \in L_T} \mu(a) y^{\dim a},
\end{equation}
where $\mu$ is the Möbius function of the lattice $L_T$.

\begin{prop}
  The characteristic polynomial of the arrangement $H_T$ is
  \begin{equation}
   \chi_T(y)= \prod_{i \in I} (y- \depth(i)),
  \end{equation}
  where $\depth(i)$ is the depth of $i$ in $T$.
\end{prop}

The number of chambers (connected components) in the real complement
is related by a theorem of Zaslavsky \cite{zaslavsky} to the value at
$y=-1$ of the characteristic polynomial .

\begin{prop}
  The number of chambers in the real arrangement $H_T$ is given by
  \begin{equation}
    \prod_{i \in I} (\depth(i)+1).
  \end{equation}
\end{prop}

\subsection{Partial order $\subseteq$ on the set of forests}

\label{subseteq}

Here is defined a partial order on the set of forests on $I$, denoted
by $\subseteq$. This poset will be shown in the next paragraphs to
contain (as intervals) all lattices $L_T$ for trees $T$ on $I$. A
similar (but different) partial order on forests has appeared in
\cite{pitman} and Exercise 5.29 of \cite{stanleybook2}.

\smallskip

Let $I$ be a finite set. A \textit{forest} on $I$ is a simply
connected graph on the vertex set $I$ where each connected component
has a distinguished vertex called its root. Therefore a forest on $I$
is a partition of $I$ together with a rooted tree on each part. The
partition of $I$ underlying a forest $F$ is denoted by $\pi(F)$. The
set of forests on $I$ is denoted by $\for(I)$. Vertices of a forest
which are not roots are called \textit{nodes}.  Let $N(F)$ be the set
of nodes of a forest $F$.

\smallskip

One can see a forest on $I$ as a structure of poset on $I$. The
partial order is the ascendance relation, \textit{i.e.} a vertex $i$
is lower than $j$ if they belong to the same tree in $F$ and $i$ is on
the path from $j$ to the root of this tree. This order relation is
denoted by $\leq_F$.  This is an extension of the definition of a tree
as a poset in \S \ref{category}.

\smallskip

Let us define the partial order $\subseteq$ on $\for(I)$ as follows.
Let $F$ and $F'$ be two forests on $I$. Then set $F' \subseteq F$ if
\begin{itemize}
\item the partition $\pi(F')$ is finer than $\pi(F)$,
\item each tree of $F'$ is induced by $F$ as a poset.
\end{itemize}
This clearly defines a partial order relation.

\begin{lemma}
  \label{inclunode}
  If $F' \subseteq F $, then $N(F')$ is contained in $N(F)$.
\end{lemma}
\begin{proof}
  A root in $F$ is a minimal element for $\leq_F$ so it is also a
  minimal element for $\leq_{F'}$ \textit{i.e.} a root in $F'$.
\end{proof}

\begin{lemma}
  \label{rankprop}
  If $F' \subseteq F $ and $N(F')=N(F)$, then $F=F'$.
\end{lemma}
\begin{proof}
  If $F$ and $F'$ have the same nodes, they have the same roots, hence
  the same number of connected components. Therefore their partitions
  are the same, hence $F=F'$.
\end{proof}

The poset $(\for(I),\subseteq)$ has a unique minimal element, which is
the forest made uniquely of roots and will be denoted by $\zero$.

\subsection{Intervals for $\subseteq$}

Let $T$ be a fixed tree on $I$. This section gives a simple
description of the interval between $\zero$ and $T$ in the poset
$(\for(I),\subseteq)$.

\begin{lemma}
  A forest $F\subseteq T$ is uniquely determined by its partition
  $\pi(F)$. 
\end{lemma}

\begin{proof}
  Given $\pi(F)$, one can reconstruct $F$ by inducing the partial
  order $\leq_T$ on the parts of $\pi(F)$.
\end{proof}

\begin{lemma}
  A partition $\pi$ can be written $\pi(F)$ for a forest $F\subseteq
  T$ if and only if each part of $\pi$ has a unique minimal element
  for the partial order induced by $\leq_T$.
\end{lemma}

\begin{proof}
  If $\pi=\pi(F)$, then each part of $\pi$ is a tree, so has a unique
  minimal element. Conversely, if a part $c$ of $\pi$ has a unique
  minimal element, then in fact it is a tree, because a poset induced
  by a tree is either a tree or a forest.
\end{proof}

\begin{lemma}
  Let $F\subseteq T$ and $F' \subseteq T$. Then $F \subseteq F'$ if
  and only if $\pi(F)$ is finer than $\pi(F')$.
\end{lemma}

\begin{proof}
  By definition of the partial order $\subseteq$, if $F \subseteq F'$
  then $\pi(F)$ is finer than $\pi(F')$. Conversely, assume that
  $\pi(F)$ is finer than $\pi(F')$. As $F$ and $F'$ are recovered from
  their partitions by inducing $\leq_T$, $F$ is in fact induced from
  $F'$, so $F\subseteq F'$.
\end{proof}

\subsection{Lattices $L_T$ and forests}

\label{descri}

In this section, one obtains a description of the lattice $L_T$ in
terms of forests.

\begin{prop}
  The elements of $L_T$ are in bijection with the forests $F$ on $I$
  which satisfy $F \subseteq T$. 
\end{prop}

\begin{proof}
  Let $a$ be an element of $L_T$, \textit{i.e.} the intersection of
  some hyperplanes $\alpha_{i,j}$ with $i \leq_T j$.
  
  Define a relation $\leq_a$ on $I$ by setting $i\leq_a j$ if
  $\alpha_{i,j}$ vanishes on $a$. Remark that $i \leq_a j$ implies $i
  \leq_T j$. That $\leq_a$ is a partial order follows from relations
  (\ref{elementaire}).  Let $\pi_a$ be the partition of $I$ given by
  the connected components of $\leq_a$.
  
  Now consider a part $c$ of $\pi_a$. Let us prove that $c$ has a
  unique minimal element. Let $i$ and $j$ be two minimal elements of
  $c$ and assume additionally that there exists $k$ in $c$ such that
  $i \leq_a k$ and $j \leq_a k$. One has either $i\leq_T j$ or $j
  \leq_T i$. Then relations (\ref{elementaire}) imply that $i \leq_a
  j$ or $j \leq_a i$, so that in fact $i=j$. Now any two minimal
  elements $i$ and $j$ in $c$ can be connected by an alternating chain
  \begin{equation}
    i \leq_a k_0 \geq_a i_0 \leq_a k_1 \geq_a i_1 \leq_a \dots \leq_a
    k_N 
    \geq_a j,
  \end{equation}
  where one can assume without restriction that $i_0,\dots,i_{N-1}$
  are minimal elements. Then a repeated application of the preceding
  argument implies that $i=i_0=i_1=\dots=i_{N-1}=j$. Therefore each
  part $c$ of $\pi_a$ has a unique minimal element for $\leq_a$. 
  
  The partial order $\leq_a$ on each part $c$ is induced by $\leq_T$.
  Indeed let $i,j$ be two elements of $c$. As said before, if $i
  \leq_a j$, then $i \leq_T j$. Conversely, assume that $i \leq_T
  j$. Let $k$ be the minimal element of $c$ for $\leq_a$. Then $k
  \leq_a i$ and $k\leq_a j$, so relations (\ref{elementaire})
  implies that $i \leq_a j$.
  
  Therefore to each element $a$ of $L_T$ is associated a forest $F_a
  \subseteq T$, which is defined as a poset by $i \leq_{F_a} j$ if and
  only if $\alpha_{i,j}$ vanishes on $a$.
  
  Conversely, one can map each forest $F \subseteq T$ to the
  intersection $a_F$ of the hyperplanes $\alpha_{i,j}$ for $i \leq_F
  j$. This set of linear forms is closed with respect to the relations
  (\ref{elementaire}). Therefore the linear forms of $H_T$ vanishing
  on $a_F$ are exactly the $\alpha_{i,j}$ for $i \leq_F j$.
  
  Furthermore, it is easy to see that the element of $L_T$ associated
  in this way to the forest $F_a$ is exactly $a$.

  This gives the sought-for bijection.
\end{proof}

\begin{theo}
  The interval between the minimal forest $\zero$ and a tree $T$ is
  isomorphic to the lattice $L_T$, \textit{i.e.} the reverse inclusion
  order in $L_T$ is mapped by the bijection to the relation
  $\subseteq$.
\end{theo}

\begin{proof}
  Let $a$ and $b$ in $L_T$ and let $F_a$ and $F_b$ be the
  corresponding forests.
  
  That $a$ is smaller than $b$ in $L_T$ means that each equation
  $\alpha_{i,j}=0$ satisfied in $a$ is also satisfied in $b$. This
  imply that the partition $\pi(F_a)$ is finer than the partition
  $\pi(F_b)$.
  
  Conversely, if the partition $\pi(F_a)$ is finer than the
  partition $\pi(F_b)$, then each tree of $a$ is induced by a tree of
  $b$ as a poset hence each relation satisfied in $a$ is also
  satisfied in $b$.
\end{proof}

\begin{figure}
  \begin{center}
    \leavevmode 
    \epsfig{file=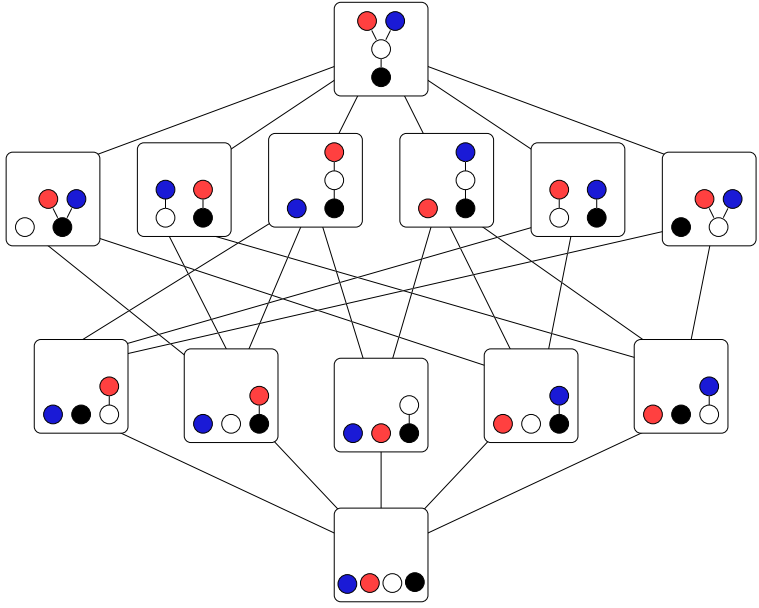,width=10cm} 
    \caption{A lattice $L_T$.}
    \label{fig-latt}
  \end{center}
\end{figure}

For an example of lattice $L_T$, see Fig. \ref{fig-latt}. 

This result shows that there is some similarity between the poset
$(\for(I),\subseteq)$ on rooted vertex-labeled forests and the poset
on forests of leaf-labeled binary trees introduced in \cite{bessel}
and further studied in \cite{besselposet}, notably because the
characteristic polynomials of intervals have a nice factorization in
both cases. This will motivate the construction of a coalgebra in \S
\ref{coalgebra}.

\begin{prop}
  The dimension of a element $a \in L_T$ is mapped by the bijection to
  the number of nodes of the associated forest $F_a$.
\end{prop}

\begin{proof}
  It is known that $L_T$ is a ranked lattice. All maximal chains have
  the same length and contain one element in each dimension. On the
  other hand, consider a maximal chain in the interval $[\zero,T]$ in
  $\for(I)$. As the set of nodes must grow at each step of the chain
  by Lemma \ref{rankprop},
  it must grow by one element only. From this follows the lemma.
\end{proof}

\begin{lemma}
  \label{treillis}
  The interval between $\zero$ and a forest $F$ is a lattice.
\end{lemma}
\begin{proof}
  There exists a tree $T$ such that $F \subseteq T$. To build one such
  $T$, choose one root $r$ of $F$ and glue the other roots of $F$ to
  some vertices of the tree with root $r$. Therefore the interval
  $[\zero,F]$ is an interval in the lattice $[\zero,T]$.
\end{proof}

\begin{lemma}
  \label{suprema}
  Let $T$ be a tree, $F_1$ and $F_2$ be forests in $L_T$ and $F$ be
  the supremum of $F_1$ and $F_2$. Assume that $N(F_1) \cap N(F_2)=
  \emptyset$. Then $N(F)=N(F_1) \sqcup N(F_2)$.
\end{lemma}

\begin{proof}
  Let $N_1=N(F_1)$ and $N_2=N(F_2)$ for short. Let $R$ be the
  complement of $N_1 \sqcup N_2$ in $I$. By Lemma \ref{inclunode}, one
  has an inclusion $ N_1 \sqcup N_2 \subseteq N(F)$.
 
  As the lattice $L_T$ is graded by the number of nodes, it is enough
  to prove that there exists a forest $F$ greater than $F_1$ and $F_2$
  with nodes $N_1 \sqcup N_2$.
  
  Consider the partition $\pi$ of $I$ which is the sup of the
  partitions $\pi_1$ and $\pi_2$ associated to $F_1$ and $F_2$
  respectively.
  
  Let $c$ be an arbitrary part of $\pi$ endowed with the partial order
  $\leq$ induced from $\leq_T$. Then any minimal element of $c$ is in
  $R$. Assume on the contrary that it is in $N_1$ for example. Then it
  cannot be minimal in $c$, for it is already not minimal in the part
  of $\pi_1$ in which it is contained, because it is a node there. The
  same is true for $N_2$ by symmetry.
  
  Now take $r \in R$. As $I$ is finite, the part $c$ of $\pi$
  containing $r$ can be built by iterated closures of $\{ r\}$ with
  respect to the partitions $\pi_1$ and $\pi_2$. Closing $\{ r\}$ with
  respect to $\pi_1$ can only add elements of $N_1$. Then closing with
  respect to $\pi_2$ can only add elements of $N_2$. This goes on in
  the same way, adding in an alternating way elements of $N_1$ and
  $N_2$ until the full part $c$ of $\pi$ containing $r$ is reached. This
  implies that each part of $\pi$ contains an unique element of $R$,
  which is its minimum.

  Therefore $\pi$ defines a forest $F$ in $L_T$, which is greater than
  $F_1$ and $F_2$ by construction and has $R$ as its set of roots. So
  $N(F)=N_1 \sqcup N_2$ and the proof is done.

\end{proof}

\subsection{Cardinality}

Here is computed the cardinality polynomial of the graded lattice
$L_T$, which counts elements of $L_T$ according to their rank. More
precisely, a recursion is found for a refinement of the cardinality
polynomial.

Define the refined generating polynomial
\begin{equation}
  C_T(y,z)=\sum_{F\in L_T} y^{\crk(F)}z^{\souche(F)},
\end{equation}
where $\souche(F)$ is one less than the cardinal of the part of
$\pi(F)$ containing the root of $T$ and $\crk(F)$ is the corank of $F$
(\textit{i.e.} one less than the number of roots). The value at $z=1$
of $C(y,z)$ is the cardinality polynomial.

Let $T$ be a tree on $I$. Let $o(T)$ be the tree on $I \sqcup \{j\}$
obtained from $T$ by grafting the root of $T$ on a new root $j$.

\begin{prop}
  \label{croissance}
  One has
  \begin{equation}
  C_{o(T)}(y,z)=z C_{T}(y,z)+y C_{T} (y,1+z).
  \end{equation}
\end{prop}
\begin{proof}
  Let $i \in I$ be the root of $T$. Let $F$ be a forest in $L_{o(T)}$.
  
  Assume first that $F$ does not have a tree with root $i$. The set of
  such forests in $L_{o(T)}$ is in bijection with the set of forests
  in $L_{T}$ in the following way. Necessarily $i$ belongs to a tree
  of $F$ with root $j$, and $i$ is the only vertex related by an edge
  to $j$. By removing $j$, one can therefore define a forest $F'$,
  with a tree of root $i$. This forest $F'$ is in $L_T$. Conversely,
  by grafting back $j$ under the root $i$ in any forest $F'$ in
  $L_{T}$, one gets an element $F$ of $L_{o(T)}$. This gives the
  sought-for bijection.
  
  Through this bijection, the number of roots is unchanged between $F$
  and $F'$. The number of nodes in the component of $j$ in $F$ is one
  more than the number of nodes in the component of $i$ in $F'$. This
  gives the term $z C_{T}(y,z)$.
  
  Assume now that $F$ does have a tree with root $i$. The set of such
  forests in $L_{o(T)}$ is in bijection with the set of pairs $(F',S)$
  where $F'$ is a forest in $L_{T}$ and $S$ is a subset of the set of
  nodes in the component of $i$ in $F'$. Consider the partition of $I$
  defined by gathering the parts containing $i$ and $j$ in the
  partition of $I \sqcup \{j\}$ defined by $F$, then removing $j$.
  Taking the order induced by $\leq_T$ on each part of this partition
  define a forest $F'$ in $L_T$. The forest $F'$ differs from $F$ only
  by the replacement of the trees $t_i$ with root $i$ and $t_j$ with
  root $j$ in $F$ by a tree $t'_i$ of root $i$ in $F'$. The set $S$ is
  defined as the set of nodes in the component of $j$ in $F$.
  Conversely, one can recover $F$ from the data of $F'$ and $S$. It is
  enough to recover the trees $t_i$ and $t_j$. The tree $t_i$ of $F$
  is induced by the tree $t'_i$ on the set of vertices not in $S$.
  The tree $t_j$ can be defined by replacing $i$ by $j$ in the tree
  induced on $\{ i\} \sqcup S$ by the tree $t'_i$.
  
  The number of roots in $F$ is one more than the number of roots in
  $F'$. The number of nodes in the component of $j$ in $F$ is the
  cardinal of $S$. This gives the term $y C_{T} (y,1+z)$.
\end{proof}

\begin{figure}
  \begin{center}
    \leavevmode 
    \epsfig{file=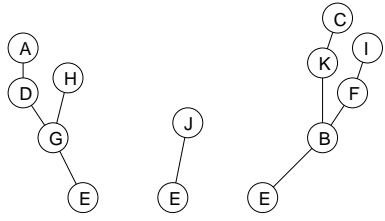,width=5cm} 
    \caption{Product decomposition.}
    \label{fig-prod}
  \end{center}
\end{figure}

On the other hand, the lattice $L_T$ is isomorphic to the product of
the lattices $L_{T_k}$ where the $T_k$ are the trees with root of
valence one obtained by separately grafting back the different
subtrees of the root. For example, Fig. \ref{fig-prod} depicts this
decomposition for the rooted tree of Fig. \ref{fig-prof}.

\begin{prop}
  \label{prop-prod}
  One has
  \begin{equation}
    C_T(y,z)=\prod_k C_{T_k}(y,z).
  \end{equation}
\end{prop}

\begin{proof}
  The bijection between $L_T$ and $\prod_k L_{T_k}$ can be described
  as follows. An element of $L_T$ is a forest $F$ with a unique tree
  containing the root of $T$. One can decompose this tree into trees
  contained in $T_k$ as in Fig. \ref{fig-prof}. Each other tree of $F$
  is contained in some unique $T_k$. Collecting trees according to the
  $T_k$ in which they are contained, one gets a collection of forests
  in $L_{T_k}$, which is the image of $F$ in $\prod_k L_{T_k}$. The
  statement of the proposition follows easily.
\end{proof}

Together, Prop. \ref{croissance} and Prop. \ref{prop-prod} give a
recursive procedure for computing the refined cardinality polynomial.

%%%%%%%%%%%%%%%%%%%%%%%%%%%%%%%%%%%%%%%%%%%%%%%%%%%%%%%%%%%%%%%%%%%%%%
                                % partie III 
%%%%%%%%%%%%%%%%%%%%%%%%%%%%%%%%%%%%%%%%%%%%%%%%%%%%%%%%%%%%%%%%%%%%%%

\section{Coalgebras on forests}

\label{coalgebra}

In this section is defined a coalgebra, using the poset $\for(I)$.

Let $F$ be a forest on $I$ and $N$ be a subset of $N(F)$. Define
\begin{equation}
  \gamma(F,N)=\sum_{F' \subseteq F \atop{N(F')=N}} F'.
\end{equation}

As this sum is without multiplicity, $\gamma(F,N)$ can also be seen as
a set.

From now on, modules, algebras and coalgebras are in the category of
$\ZZ$-graded abelian group with Koszul sign rules for the tensor
product.

An orientation $o$ of a forest $F$ is a maximal exterior product of
the set $N(F)\sqcup \{R_F\}$, where $R_F$ is an auxiliary element. An
oriented forest is a tensor product $o\otimes F$, where $o$ is an
orientation of the forest $F$.

\smallskip

Let $\For(I)$ be the free $\ZZ$-module on the set of oriented forests
on $I$ modulo the relations $(-o)\otimes F=-(o \otimes F)$.  This
module is graded by $\deg(o \otimes F)=\# N(F)$. A map $\Delta$ from
$\For(I)$ to $\For(I) \otimes \For(I)$ is defined as follows:
\begin{equation}
  \label{deficop}
  \Delta(o\otimes F)=
  \sum_{N(F)=N_1 \sqcup N_2} (o_1 \otimes \gamma(F,N_1))
  \otimes(o_2 \otimes \gamma(F,N_2)),
\end{equation}
where the orientations satisfy $o=o_1 \wedge r \wedge o_2$ modulo
$R=R_1 \wedge r \wedge R_2$.

\begin{lemma}
  \label{cruxasso}
  Let $F$ be a forest in $L_T$. Let $N_1$ be a subset of $N(F)$ and
  fix a partition $N_1=N'_1 \sqcup N''_1$. Then (not taking care of
  orientations), one has
  \begin{equation}
    \sum_{F_1 \in \gamma(F,N_1)} \gamma(F_1,N'_1)\otimes \gamma(F_1,N''_1) 
    = \gamma(F,N'_1)\otimes \gamma(F,N''_1).
  \end{equation}
\end{lemma}

\begin{proof}
  It is sufficient to give a bijection between the
  set of triples $(F'_1,F''_1,F_1)$ satisfying
  \begin{align}
    F_1 &\in \gamma(F,N_1),\\
    F'_1 &\in \gamma(F_1,N'_1),\\
    F''_1 &\in \gamma(F_1,N''_1),
  \end{align}
  and the set of pairs $(F'_1,F''_1)$ satisfying
  \begin{align}
    F'_1 &\in \gamma(F,N'_1),\\
    F''_1 &\in \gamma(F,N''_1).
  \end{align}

  In one direction, a triple is mapped to a pair by forgetting $F_1$.
  
  In the other direction, let $F'_1$ and $F''_1$ be a pair as above.
  By Lemma \ref{suprema}, the supremum $F_1$ of $F'_1$ and $F''_1$
  (well-defined by Lemma \ref{treillis}) satisfies the conditions
  $F'_1 \subseteq F_1$, $F''_1 \subseteq F_1$ and $N(F)=N_1$. By Lemma
  \ref{rankprop}, this is the unique forest satisfying these
  conditions.
\end{proof}

\begin{prop}
  Formula (\ref{deficop}) endows $\For(I)$ with the structure of a
  cocommutative coassociative counital coalgebra with coproduct
  $\Delta$.
\end{prop}

\begin{proof}
  Cocommutativity is clear from Formula (\ref{deficop}).
  Coassociativity is deduced from the following formula for the double
  coproduct, which is a consequence of Lemma \ref{cruxasso}: (details
  are left to the reader)
  \begin{equation}
    \sum_{N_F=N_1 \sqcup N_2 \sqcup N_3} (o_1 \otimes \gamma(F,N_1))
    \otimes(o_2 \otimes \gamma(F,N_2))\otimes(o_3 \otimes \gamma(F,N_3)),
  \end{equation}
  where $o=o_1 \wedge r \wedge o_2 \wedge s \wedge o_3$ modulo $R=R_1
  \wedge r \wedge R_2 \wedge s \wedge R_3$. The counit is the
  projection to the degree-zero component, which is one-dimensional.
\end{proof}

\begin{lemma}
  \label{disnode}
  Let $F$ be a forest on $I$ and $F_1 \otimes F_2$ be a term in the
  coproduct $\Delta(F)$. Then the set of nodes of $F_1$ and $F_2$ are
  disjoint in $I$.
\end{lemma}

\begin{proof}
  Obvious from Formula (\ref{deficop}).
\end{proof}

\begin{lemma}
  \label{ascend}
  Let $F$ be a forest on $I$ and $F_1 \otimes F_2 \otimes \dots
  \otimes F_k$
  be a term in the iterated coproduct of $F$. Assume that one has $i
  \leq_{F_\ell} j$ for some $\ell$. Then one has $i\leq_F j $. In
  words, each ascendance relation in a term of a coproduct of $F$ is
  also satisfied in $F$.
\end{lemma}

\begin{proof}
  By definition of the partial order $\subseteq$, each tree in each
  $F_{\ell}$ is induced by $F$ as a poset.
\end{proof}

\section{Algebras on forests}

\subsection{Partial order $\preceq$ on forests}

Consider the following partial order $\preceq$ on the set of forests
on $I$. One sets $F \preceq F'$ if the identity of $I$ is a morphism
of posets from $(I,\leq_F)$ to $(I,\leq_{F'})$. Note that in order to
check if $F \preceq F'$, it is enough to check the relations $i
\leq_{F'} j$ for all edges $i \leftarrow j$ of $F$ .

Observe that this partial order is completely different from the
partial order $\subseteq$ introduced in \S \ref{subseteq}. For
example, there are distinct forests with the same nodes which are
comparable for $\preceq$ (compare to Lemma \ref{rankprop}).

\subsection{Generators and relations}

Recall that the ambient category is that of $\ZZ$-graded abelian
groups with Koszul sign rules. So there are appropriate signs to be
inserted whenever two elements are exchanged.

In this category, consider the commutative associative unital algebra
$\mon(I)$ defined by generators $\Omega_{i,j}$ of degree $1$ for
$i\not=j$ in $I$ and relations
\begin{equation}
  \label{fourche}
  \Omega_{i,k} \Omega_{j,k}=0,
\end{equation}
for $i,j,k$ pairwise different in $I$, and
\begin{equation}
  \label{boucle}
  \Omega_{i_0,i_1} \Omega_{i_1,i_2} \dots \Omega_{i_n,i_0}=0,
\end{equation}
for $n \geq 1$ and pairwise different $i_0,\dots,i_n$ in $I$.

\begin{prop}
  The algebra $\mon(I)$ is spanned by monomials $m_F$ indexed by the
  set $\for(I)$ of forests of rooted trees on $I$.
\end{prop}

\begin{proof}
  One can represent (up to sign) a monomial in the generators as an
  oriented graph on $I$ by drawing an oriented edge $i \leftarrow j$
  for the generator $\Omega_{i,j}$.
  
  The relation (\ref{fourche}) says exactly that graphs containing two
  edges going out of the same vertex vanish. Therefore any cycle in a
  non-vanishing graph is oriented.
  
  The relations (\ref{boucle}) say exactly that graphs containing an
  oriented cycle vanish.
  
  It follows that the algebra $\mon(I)$ is spanned by monomials
  indexed by oriented simple graphs on $I$ with no cycle and no
  divergence of arrows, \textit{i.e.} rooted forests on $I$. For each
  forest $F$, choose (well defined up to sign) a monomial $m_F$. 
\end{proof}

In fact, this set of monomials $(m_F)_F$ form a basis of $\mon(I)$.
This property can be seen directly or will follow from Theorem
\ref{isotheo}.

\subsection{Isomorphism}

Consider now the dual algebra $\For^*(I)$ of the coalgebra $\For(I)$,
for the pairing $\langle \, \rangle : \For^*(I)\otimes \For(I) \to
\ZZ$. Let $(F^*)_F$ be the basis of $\For^*(I)$ dual to the basis
$F$ of $\For(I)$ for this pairing. 

Define elements $F^*_{i,j}$ of degree $1$ in $\For^*(I)$ by
\begin{equation}
  \langle F^*_{i,j},F_{k,\ell} \rangle = \delta_{i,k}\delta_{j,\ell},
\end{equation}
where $F_{k,\ell}$ is the unique forest of degree $1$ on $I$ with only
one node, labeled $\ell$ and attached to a root labeled $k$. The
orientation of $F_{k,\ell}$ is prescribed by $\ell \wedge R$.

\begin{theo}
  \label{isotheo}
  The algebra $\For^*(I)$ is isomorphic to the algebra $\mon(I)$.
\end{theo}

First, one has to define a map from $\mon(I)$ to $\For^*(I)$.

\begin{prop}
  There is a map $\rho$ from $\mon(I)$ to $\For^*(I)$ satisfying
  \begin{equation}
    \rho( \Omega_{i,j} ) = F^*_{i,j}.
  \end{equation}
\end{prop}

\begin{proof}
  One has to show that the elements $F^*_{i,j}$ of $\For^*(I)$ satisfy the
  relations (\ref{fourche}) and (\ref{boucle}).
  
  To check the relation (\ref{fourche}), it is sufficient to prove
  that the term $F_{i,k} \otimes F_{j,k}$ does not appear in the
  coproduct $\Delta(F)$ of any forest $F$ on $I$. This follows from
  Lemma \ref{disnode}, because $k$ is a node of both $F_{i,k}$ and
  $F_{j,k}$.

  To check the relation (\ref{boucle}), it is sufficient to prove
  that the term 
  \begin{equation}
    F_{i_0,i_1} \otimes F_{i_1,i_2} \otimes \dots \otimes F_{i_n,i_0} 
  \end{equation}
  does not appear in the iterated coproduct of any forest $F$ on $I$.
  Assume that it does appear in the iterated coproduct of $F$. Then by
  Lemma \ref{ascend}, one has
  \begin{equation}
     i_0 \leq_F i_1 \leq_F i_2 \leq_F \dots \leq_F i_n \leq_F i_0, 
  \end{equation}
  which would mean that there is a cycle in $F$, which is absurd.
  
  Therefore the map $\rho$ is well-defined.
\end{proof}

The crucial part of the proof of Theorem \ref{isotheo} is the
following Lemma.

\begin{lemma}
  \label{imagemono}
  Let $F$ be a forest on $I$. The image by $\rho$ of the monomial
  $m_F$ is the sum (with coefficients $\pm 1$) of $G^*$ over the set
  of all forests $G$ which are $\succeq$ than $F$ and have the same
  nodes as $F$.
\end{lemma}

\begin{proof}
  Write $m_F$ as a product $\Omega_{i_0,j_0} \dots \Omega_{i_n,j_n}$
  over the set of edges of $F$, in some order. Then $\rho(m_F)$ is
  $F^*_{i_0,j_0} \dots F^*_{i_n,j_n}$.
  
  Let $G$ be a forest on $I$. The coefficient of $G^*$ in
  $F^*_{i_0,j_0} \dots F^*_{i_n,j_n}$ is (up to sign) the coefficient
  of $F_{i_0,j_0} \otimes \dots \otimes F_{i_n,j_n}$ in the iterated
  coproduct of $G$.
  
  Assume that the coefficient of $G^*$ in $\rho(m_F)$ is non-zero.
  Then the coefficient of $F_{i_0,j_0} \otimes \dots \otimes
  F_{i_n,j_n}$ in the iterated coproduct of $G$ is non-zero. Lemma
  \ref{ascend} implies, for each edge $i \leftarrow j$ of $F$, that $i
  \leq_{G} j$. Hence by definition of $\preceq$, one has $F \preceq
  G$.
  
  By homogeneity of the product, one also has that $G$ and $F$ have
  the same number of nodes and the same number of roots. As $F \preceq
  G$, each root of $F$ is also a root of $G$, so $F$ and $G$ have the
  same roots and nodes.
  
  Now consider any forest $G$ with the same nodes and roots as $F$ and
  satisfying $F \preceq G$. For such a $G$, the tensor $F_{i_0,j_0}
  \otimes \dots \otimes F_{i_n,j_n}$ does appear with coefficient $\pm
  1$ in the iterated coproduct of $G$, in the term of the sum
  (\ref{deficop}) corresponding to the partition of
  $N(G)=\{j_0,j_1,\dots,j_n\}$ into singletons. Indeed, for any $k$,
  $F \preceq G$ implies $i_k \leq_G j_k$, which in turn implies that
  $F_{i_k,j_k}$ belongs to the set $\gamma(G,\{j_k\})$.
\end{proof}

The proof of Theorem \ref{isotheo} can now be completed.

\begin{proof}
  By Lemma \ref{imagemono}, the image by $\rho$ of the set $(m_F)_F$
  which spans $\mon(I)$ is a basis of $\For^*(I)$. As the rank of
  $\For^*(I)$ is the cardinal of $\for(I)$, one deduces that $(m_F)_F$
  is in fact a basis of $\mon(I)$. So $\rho$ is an isomorphism.
\end{proof}

%%%%%%%%%%%%%%%%%%%%%%%%%%%%%%%%%%%%%%%%%%%%%%%%%%%%%%%%%%%%%%%%%%%%%%
\bibliographystyle{plain}
\bibliography{arbres}

\end{document}